\documentclass[11pt,amssymb]{article}
\usepackage{eucal}
\usepackage{amsfonts,amsmath,amssymb,latexsym}
\usepackage{geometry}

\begin{document}

\newtheorem{theorem}{Theorem}
\newtheorem{proposition}{Proposition}
\newtheorem{corollary}{Corollary}
\newtheorem{lemma}{Lemma}
\newtheorem{definition}{Definition}
\newtheorem{remark}{Remark}
\newtheorem{claim}{Claim}
\numberwithin{equation}{section}

\def\A{\mathcal{A}}
\def\Ss{\mathbb{S}}
\def\Mm{\mathbb{M}^{n}({\kappa})}
\def\Mmn{\mathbb{M}^{n}({\kappa})}
\def\Cc{\mathcal{C}}
\def\Rr{\mathbb{R}}
\def\Pp{\mathbb{P}_t}
\def\dd{d}
\def\mt{\mathcal}
\def\ae{\langle}
\def\ad{\rangle}
\def\sn{\textrm{sn}}
\def\ct{\textrm{ct}}
\def\cs{\textrm{cs}}
\def\re{\textrm{Re}}
\def\im{\textrm{Im}}
\def\Qq{\mathbb{Q}}

\def\esq{\langle}
\def\dir{\rangle}

\title{Existence of isometric immersions into nilpotent Lie groups}
\author{Jorge H. S. de Lira\thanks{partially supported by CNPq and FUNCAP} \,\,\,and\,\, Marcos F. Melo
\thanks{partially supported by CAPES and CNPq}}
\date{}
\maketitle

\begin{abstract}
We establish necessary and sufficient conditions for existence of
isometric immersions of a simply connected Riemannian manifold into
a two-step nilpotent Lie group. This comprises the case of
immersions into $H$-type groups.
\end{abstract}

\vspace{3mm}


\vspace{3mm}

\section{Introduction}

The fundamental theorem of submanifold theory, usually referred to
as Bonnet's theorem, states that the Gauss, Codazzi and Ricci
equations constitute a set of integrability conditions for isometric
immersions of a simply connected Riemannian manifold in Euclidean
space with prescribed second fundamental form. From the viewpoint of
exterior differential systems, this result is a classical
application of Frobenius's theorem. At this respect, we refer the
reader to \cite{dajczer}, \cite{landsberg} and \cite{spivak} for
instance.

Versions of Bonnet's Theorem for immersions in Riemannian spaces
were recently achieved by Benoit Daniel in \cite{daniel-1} and
\cite{daniel-2} and by  P. Piccione and V. Tausk in \cite{piccione}
and \cite{piccione2}. In \cite{daniel-2}, Daniel consider immersions
in three dimensional homogeneous spaces with four dimensional
isometry group as Heisenberg spaces and Berger spheres. These
ambients are regarded there as total spaces of Riemannian
submersions over constant curvature surfaces, fibered by flow lines
of a vertical Killing vector field $\xi$. It is proved that the
ambient curvature tensor expressed in terms of a frame adapted to
the immersion may be completely determined by the first and second
fundamental forms and by the normal component $\nu$ and tangencial
projection $T$ of $\xi$. Since Gauss and Codazzi equations involve
these projections, it is necessary to consider two additional first
order differential equations in $\nu$ and $T$. The augmented set of
equations is then a complete integrability condition.

In \cite{piccione2},  Piccione and Tausk prove a general existence
result for affine immersions into affine manifolds endowed with a
$G$-structure. The immersions should preserve the $G$-structure and
the ambient spaces are required to be infinitesimally homogeneous.
Roughly speaking, this last condition assures that the ambient
curvature is constant when computed in terms of frames belonging to
a $G$-reduction of the frame bundle. This method encompasses all
classical results as well as Daniel's results. Another applications
of this technique in the context of Lie groups and Lorentzian Geometry 
may be found in  \cite{Sinue}, \cite{LMa}, \cite{LP}  and
\cite{Manfio}.

The Heisenberg spaces studied in \cite{daniel-2} are nilpotent Lie
groups. Indeed, two-step nilpotent Lie groups form a distinguished
class of geometric objects which include real, complex and
quaternionic Heisenberg spaces and more generally H-type groups
(see, e.g., \cite{Kaplan},  \cite{Eberlein}, \cite{Berndt} and
\cite{Berndt2}).  These groups have some remarkable analytical
properties and appear in distinct areas as Harmonic Analysis (v.
\cite{Damek}) and General Relativity (v. \cite{Hervik}).

These remarks motivate us to raise the question of extending
Bonnet's theorem from the classical case, which corresponds to
Abelian groups, to two-step nilpotent Lie groups. Theorem \ref{main}
below yields such an extension in the spirit of the results in
\cite{daniel-2}.

A brief outline of this paper may be given as follows. Let $N$ be a
$(n+n')$-dimensional two-step nilpotent Lie group, where $n'$ is the
dimension of the center $\mathfrak{z}$ in its Lie algebra. As occurs
in \cite{daniel-2}, $\mathfrak{z}$ is spanned by left-invariant $n'$
Killing vector fields $E_{n+k}$, $k=1,\ldots, n'$, whose covariant
derivatives determine certain skew-symmetric tensors $J_k$,
$k=1,\ldots, n'$. The Section 2 is devoted to show that the
curvature tensor in $N$ may be computed in an arbitrary frame
$\{e_a\}_{a=1}^{n+n'}$ solely in terms of the tensors $J_k$ and the
projections $U^k_a =\langle E_{n+k},e_a\rangle$. The curvature form
relative to the frame $\{e_a\}_{a=1}^{n+n'}$ is given by the tensor
$Q$ defined in Section \ref{Qcurv}.

In the particular case of a frame adapted to an isometric immersion,
this implies that Gauss, Codazzi and Ricci equations are completely
written only in terms of the first and second fundamental forms and
the normal and tangential projections of the Killing vector fields
$E_{n+1},\ldots, E_{n+n'}$ and their covariant derivatives. This is
the content of Section \ref{imm-nec}.

In Sections \ref{ex-frame} and \ref{ex-imm}, we establish sufficient
conditions for immersing isometrically a simply connected Riemannian
manifold $M$ into $N$, with prescribed second fundamental form.  For
this, we consider a real Riemannian vector bundle $\mathcal{E}$ over
$M$ with rank $m'=n+n'-m$ so that the Whitney sum
$\mathcal{S}=TM\oplus\mathcal{E}$ is a trivial bundle. We define an
orthonormal global frame $\{\hat E_{a}\}_{a=1}^{n+n'}$ in
$\mathcal{S}$ and then transplante the definitions of the tensors
$J_k$ and $Q$ to this setting. This may be done in last analysis
because these tensors depend on the structural constants of $N$. We
then prove

\begin{theorem}
\label{main} a) Let $M^m$ be a Riemannian simply connected manifold
and let $\mathcal{E}$ be a real Riemannian vector bundle with rank
$m'$ so that $\mathcal{S}=TM\oplus\mathcal{E}$ is a trivial vector
bundle. Let $\hat\nabla$ and $\hat R$ be respectively the compatible
connection and curvature tensor in $\mathcal{S}$ and $\nabla$ and
$\nabla^{\mathcal{E}}$ the compatible connections induced in $TM$
and $\mathcal{E}$, respectively. We fix a global orthonormal frame
$\{\hat E_k\}_{k=1}^{n+n'}$  in $\mathcal{S}$. Define $\hat J_k$ and
$\hat Q$ as in (\ref{Jk})  and (\ref{QA}), respectively. Assume that
these fields satisfy the Gauss, Codazzi and  Ricci equations
\begin{equation}
\label{main-hyp-1} \hat R = \hat Q
\end{equation}
and the additional equations
\begin{equation}
\label{main-hyp-2}  \hat\nabla \hat E_{n+k}=-1/2 \,\,\hat J_k, \quad
k=1,\ldots, n'.
\end{equation}
Thus, there exists an isometric immersion $f:M\to N$ covered by  a
bundle isomorphism $f_*^\perp: \mathcal{E}\to TM^{\perp}_f$, where
$TM^\perp_f$ is the normal bundle along $f$ so that $f^\perp_*$ is
an isometry when restrited to the fibers and satisfies
\begin{eqnarray}
& & f^\perp_* \nabla^{\mathcal{E}}_X V = \nabla_{X}^\perp f^\perp_*
V,\quad X\in \Gamma(TM), \,\, V\in \Gamma(\mathcal{E}),\\
& & f^\perp_* II(X,Y)=\bar\nabla_{f_*X} f_*Y-f_*(\nabla_X Y), \quad
X,Y\in \Gamma(TM),
\end{eqnarray}
where $\bar\nabla$ and $\nabla^\perp$ denote, respectively, the
connections in $N$ and $TM^\perp_f$ and the tensor $II\in
\Gamma(T^*M\otimes T^*M\otimes \mathcal{E})$ is defined by
\begin{equation}
\hat\nabla_X Y = \nabla_X Y + II(X,Y),\quad X,Y\in \Gamma(TM).
\end{equation}

\vspace{2mm}

\noindent b) Let $f,\tilde f$ be two isometric immersions from $M$
to $N$ 
with second fundamental forms $II_f$ and $II_{\tilde f}$ satisfying
\begin{equation}
\label{IIi} II_f (X,Y)= \Phi II_{\tilde f}(X,Y), \quad X,Y\in
\Gamma(TM),
\end{equation}
and normal connections $\nabla^\perp$ and $\tilde\nabla^\perp$ on
the respective normal bundles $TM^\perp_f$ and $TM^\perp_{\tilde f}$
related by
\begin{equation}
\label{perpi} \Phi\nabla^{\perp}_{X}V =\tilde\nabla^\perp_X \Phi(V),
\quad V\in \Gamma(TM^\perp_f),
\end{equation}
where $\Phi:TM^\perp_f\to TM^\perp_{\tilde f}$ is a vector bundle
isomorphism satisfying
\begin{equation}
\langle \Phi(V),\Phi(W)\rangle = \langle V,W\rangle, \quad V,W\in
\Gamma(TM^\perp_f).
\end{equation}
Fixed a left-invariant frame $\{E_k\}_{k=1}^{n+n'}$ in $N$ we assume
that
\begin{equation}
\label{eq-add} \langle f_* X, E_{n+k}\rangle = \langle \tilde f_* X,
E_{n+k}\rangle,\quad X\in \Gamma(TM)
\end{equation}
and that
\begin{equation}
\label{eq-add-2} \langle V, E_{n+k}\rangle =\langle
\Phi(V),E_{n+k}\rangle, \quad  V\in \Gamma(TM^\perp_f).
\end{equation}
for $k=1,\ldots,n'$.

Then, there exists an isometry $L:N\to N$ such that $\tilde f =
L\circ f$.
\end{theorem}




The ultimate reason for refer to (\ref{main-hyp-1}) as Gauss,
Codazzi and Ricci equations is that the tensor $\hat Q$ imitates the
curvature tensor in $N$ when written in terms of a frame adapted to
an isometric immersion. We point out that imposing that
$\mathcal{S}$ is trivial allows us to give an intrinsic meaning to
the tensors $\hat J_k$. Hypothesis (\ref{main-hyp-1}) and
(\ref{main-hyp-2}) play here the same role as the construction of a
flat bundle endowed with parallel sections in the proof of the
classical case.

Our method keeps some resemblance with the proof of Bonnet's theorem
given by P. Ciarlet and F. Larsonneur in \cite{Ciarlet}. Indeed,
Theorem \ref{main} may be regarded as establishing sufficient
conditions for immersing an open set of the Euclidean space into a
two-step nilpotent Lie group.

\section{Two-step nilpotent Lie groups}\label{basic}

Let $N$ be  a Lie group  with Lie algebra $\mathfrak{n}$ and
Maurer-Cartan form $\omega_{\mathfrak{n}}$. The Levi-Civit\`a
connection of a given left-invariant metric $\ae \cdot,\cdot \ad$ on
$N$ is
\begin{equation}
\label{connection_Lie} 2\bar\nabla_E F = [E,F]-\textrm{ad}^*_E \cdot
F -\textrm{ad}^*_F\cdot E,
\end{equation}
where $E,F$ are left-invariant vector fields in $\mathfrak{n}$ and
\[
\ae \textrm{ad}^*_E\cdot F,G\ad =  \ae F,[E,G]\ad, \quad
E,F,G\in\mathfrak{n}.
\]
We suppose that $\mathfrak{n}$ may be decomposed as
$\mathfrak{n}=\mathfrak{z}\oplus\mathfrak{v}$ with
\begin{equation}
\label{mh} [\mathfrak{v},\mathfrak{v}]\subset \mathfrak{z},\quad
[\mathfrak{z},\mathfrak{n}]= \{0\},
\end{equation}
what implies that $N$ is a two-step nilpotent Lie group. Let us
denote by $n$ and $n'$ the dimensions of $\mathfrak{v}$ and
$\mathfrak{z}$, respectively. We suppose that the direct sum
$\mathfrak{n}=\mathfrak{z}\oplus\mathfrak{v}$ is orthogonal. The
relations (\ref{mh}) then yield
\begin{eqnarray}
& & \label{nabla-mm}\bar\nabla_E F =\frac{1}{2}[E,F],\quad E,F\in\mathfrak{v},\\
& & \label{nabla-mh}\bar \nabla_E Z =\bar\nabla_Z E = -\frac{1}{2} J_Z E,\quad E\in\mathfrak{v},\,\,Z\in\mathfrak{z},\\
& &\label{nabla-hh} \bar \nabla_Z Z' =0, \quad Z,Z'\in\mathfrak{z},
\end{eqnarray}
where the operator $J_Z:\mathfrak{v}\to\mathfrak{v}$ associated to a
vector field $Z\in\mathfrak{z}$ is defined by $J_Z =
\textrm{ad}^*Z$. This operator may be extended to the whole algebra
$\mathfrak{n}$ as
\begin{eqnarray}
\label{J-ext} J_Z :=-2 \bar\nabla Z.
\end{eqnarray}
It is useful to consider also the $(0,2)$ tensor field equally
denoted by $J_Z$ and defined by $J_Z(E,F)=\ae J_ZE,F\ad$.

\subsection{Some auxiliary tensors}\label{LQ}

According to the decomposition
$\mathfrak{n}=\mathfrak{v}\oplus\mathfrak{z}$, we choose an
orthonormal left-invariant frame field
\begin{equation}
\label{basis} E_1,\ldots,E_n,E_{n+1},\ldots, E_{n+n'},
\end{equation}
so that the first $n$ vectors are in $\mathfrak{v}$ and the next
$n'$ ones are in $\mathfrak{z}$. Fixed this choice of frame,  we
define the structural constants of $N$ by
\begin{equation}
\label{sigma0} [E_k, E_l]=\sum_{r=1}^{n+n'}\sigma^r_{kl}E_r, \quad
1\le k,l\le n+n'.
\end{equation}
If $\{\theta^k\}_{k=1}^{n+n'}$ denotes  the co-frame dual to
$\{E_k\}_{k=1}^{n+n'}$, then the corresponding connection forms in
$N$ are given by
\begin{equation}
\label{connection_Lie_2} \theta^k_l
=\frac{1}{2}\sum_{r=1}^{n+n'}\tau^k_{lr}\, \theta^r,
\end{equation}
where
\begin{equation}
\label{tau}
 \tau^k_{lr}=\sigma^k_{rl}
+\sigma_{kr}^l + \sigma_{kl}^r.
\end{equation}
We also define the curvature $2$-form
$\Theta=\{\Theta^k_l\}_{k,l=1}^{n+n'}$ of $N$ associated to
(\ref{basis}) by
\begin{equation}
\label{curv-amb} \Theta^k_l
=\frac{1}{4}\sum_{r,s,t=1}^{n+n'}\big(\tau^k_{lr}\tau^r_{st}+\tau^k_{rs}\tau^r_{lt}\big)\theta^s\wedge\theta^t.
\end{equation}
These forms satisfy the structural equations
\begin{equation}
\label{first-N} \dd \theta^k + \sum_{l=1}^{n+n'}\theta^k_l \wedge
\theta^l =0, \quad \theta^k_l=-\theta^l_k
\end{equation}
and
\begin{equation}
\label{G-curv} \dd \theta^k_l + \sum_{r=1}^{n+n'}\theta^k_r \wedge
\theta^r_l = \Theta^k_l,
\end{equation}
where $1\le k,l\le n+n'$.

\subsubsection{Christoffel tensor}\label{ctensor}

Given the left-invariant frame above, we denote $J_k=J_{E_{n+k}}$,
$1\le k\le n'$. Fixed this notation,  we define in $N$ the tensor
field
\begin{eqnarray}
\label{L} & & L(X,Y,V) = -\frac{1}{2}\sum_{k=1}^{n'}\ae J_{k} V,X\ad
\ae Y,E_{n+k}\ad+ \frac{1}{2}\sum_{k=1}^{n'}\ae J_{k} Y,X\ad \ae
V,E_{n+k}\ad \nonumber \\
& & \,\,\,\,\,\,-\frac{1}{2}\sum_{k=1}^{n'}\ae J_{k} Y,V\ad \ae
X,E_{n+k}\ad,\quad X, Y, V\in \Gamma(TN).
\end{eqnarray}
In order to derive a {\it local} expression for $L$,  we consider a
frame $\{e_a\}_{a=1}^{n+n'}$ defined in an open set $N'$ of $N$  by
\begin{equation}
\label{eE} e_a = \sum_{b=1}^{n+n'}E_b \,A^b_a,
\end{equation}
for some map $A:N'\to \textsc{SO}_{n+n'}$. For $1\le a \le n+n'$ and
$1\le k\le n'$, we define the functions
\begin{equation}
\label{U} U^{k}_a = \theta^{n+k}(e_a) = A_a^{n+k}.
\end{equation}
Thus, if $(\omega^a)_{a=1}^{n+n'}$ and
$(\omega_{a}^b)_{a,b=1}^{n+n'}$ are respectively the dual forms and
the connection forms associated to the frame $\{e_a\}_{a=1}^{n+n'}$,
one has
\begin{eqnarray}
\label{u} \omega^a(\bar\nabla E_{n+k}) = \dd U^k_a-\sum_c
U^k_c\omega^c_a=:\frac{1}{2}\sum_b u^k_{ab}\,\omega^b.
\end{eqnarray}
Hence,  one gets
\begin{equation}
J_{k}= \sum_{a,b=1}^{n+n'} u^k_{ab}\omega^a \otimes \omega^b.
\end{equation}
Notice that
\begin{eqnarray}
\ae J_k V,W\ad & = & -2\ae \bar\nabla_{V} E_{n+k},W\ad =
-2\sum_{l,r}\ae V,E_l\ad\ae W,E_r\ad \ae \bar\nabla_{E_l}E_{n+k},E_r\ad \nonumber\\
&=& \sum_{l,r}\ae V,E_l\ad \ae W,E_r\ad
\sigma^{n+k}_{lr}\label{Jk0}.
\end{eqnarray}
In local terms, that is, setting $V=e_a, W=e_b$, one has
\begin{equation}
\label{u-2} u^k_{ab}=\sum_{l,r=1}^n A^l_a A^r_b \sigma^{n+k}_{lr}.
\end{equation}
Using the local frame, one computes
\begin{eqnarray*}
L(X,e_a,e_b) & = & -\frac{1}{2}\sum_{k=1}^{n'}\ae J_{k} e_b,X\ad \ae
e_a,E_{n+k}\ad+\frac{1}{2}\sum_{k=1}^{n'}\ae J_{k}e_a,X\ad \ae
e_b,E_{n+k}\ad
\\
& &\,\,-\frac{1}{2}\sum_{k=1}^{n'}\ae J_{k}e_a,e_b\ad \ae X,E_{n+k}\ad\\
& = &
-\frac{1}{2}\sum_{c=1}^{n+n'}\sum_{k=1}^{n'}\big(U^k_au^k_{bc}-U^k_b
u^k_{ac}+U^k_c u^k_{ab}\big)\omega^c(X).
\end{eqnarray*}
One then  defines the matrix of $1$-forms $\lambda=
(\lambda^a_b)_{a,b=1}^{n+n'}$ by
\begin{equation}
\label{lambda} \lambda^a_b = L(\,\cdot\, , e_a, e_b),
\end{equation}
that is,
\begin{equation}
\label{u3} \lambda^a_b =-\frac{1}{2}\sum_{c=1}^{n+n'}\sum_{k=1}^{n'}
\big(U^k_a u^k_{bc}-U^k_bu^k_{ac}+U^k_c u^k_{ba}\big)\omega^c.
\end{equation}

We now use the equation (\ref{u-2}) for obtaining an alternative
expression for $\lambda$, which will be useful later.

\begin{proposition}\label{l}
The $1$-form $\lambda=(\lambda^a_b)_{a,b=1}^{n+n'}$ defined in
(\ref{lambda}) satisfy
\begin{equation}
 \lambda=A^{-1}\theta A,
\end{equation}
where $\theta=(\theta^k_l)_{k,l=1}^{n+n'}$. Thus, the connection
form $\omega=(\omega^a_b)_{a,b=1}^{n+n'}$ is given by
\begin{equation}
\omega = A^{-1}\dd A +\lambda.
\end{equation}
\end{proposition}

\noindent{\it Proof.} Using (\ref{connection_Lie}) and (\ref{tau}),
one obtains
$$
\tau_{lr}^k=\left\{
\begin{array}{l}
\sigma_{rl}^k,\ 1\leq l,r\leq n \ \textrm{and} \ k\geq n+1,\\
\sigma_{kl}^r,\ 1\leq k,l\leq n \ \textrm{and} \ r\geq n+1,\\
\sigma_{kr}^l,\ 1\leq k,r\leq n \ \textrm{and} \ l\geq n+1.\\
\end{array}
\right.
$$
Thus, (\ref{u3}) and (\ref{u-2}) yield
\begin{eqnarray*}
\lambda_b^a
& =&-\frac{1}{2}\sum_{c=1}^{n+n'}\sum_{k=1}^{n'}(U_a^ku_{bc}^k-U_b^ku_{ac}^k-U_c^ku_{ab}^k)\omega^c\\
&   = &
-\frac{1}{2}\sum_{c=1}^{n+n'}\sum_{k=1}^{n'}\sum_{l,r=1}^nU_a^kA_b^lA_c^r\sigma_{lr}^{n+k}\omega^c+
\frac{1}{2}\sum_{c=1}^{n+n'}\sum_{l=1}^{n'}\sum_{k,r=1}^nA_a^kU_b^lA_c^r\sigma_{kr}^{n+l}\omega^c\\
& &\,
+\frac{1}{2}\sum_{c=1}^{n+n'}\sum_{r=1}^{n'}\sum_{k,l=1}^nA_a^kA_b^lU_c^r\sigma_{kl}^{n+r}\omega^c
\end{eqnarray*}
what implies that
\begin{eqnarray*}
\lambda_b^a &  = &\frac{1}{2}\sum_{c,k,l,r=1}^{n+n'}
A_a^kA_b^lA_c^r\tau_{lr}^k\omega^c
=\frac{1}{2}\sum_{k,l,r=1}^{n+n'}A_a^kA_b^l\tau_{lr}^k\theta^r
=\sum_{k,l=1}^{n+n'}A_a^kA_b^l\theta_l^k\\
&  = &
\sum_{k,l=1}^{n+n'}(A^{-1})_k^a\theta_l^kA_b^l\\
&=&(A^{-1}\theta A)_b^a
\end{eqnarray*}
This finishes the proof of the proposition. \hfill $\square$

\subsubsection{A curvature-type tensor}\label{Qcurv}

We then define a $(0,4)$ covariant tensor $Q$ in $N$ by
\begin{equation}
\label{Qamb} Q(X,Y,V,W)=Q_1(X,Y,V,W)+Q_2(X,Y,V,W),
\end{equation}
where $X,Y,V$ and $W$ are vector fields in $N$ and $Q_1$ and $Q_2$
are given by
\begin{eqnarray*}
& & Q_1(X, Y, V, W)\\
& &\,\, =\frac{1}{4}\ae J_k X,W\ad\ae J_k V,Y\ad +\frac{1}{2}\ae J_k
Y,X\ad \ae J_k W,V\ad-\frac{1}{4}\ae J_k Y,W\ad \ae
J_k V,X\ad\\
& & \,\,\,\,\,-\frac{1}{2}\sum_k\ae W,E_{n+k}\ad \ae(\bar\nabla_X
J_k)V,Y\ad
+\frac{1}{2}\sum_k\ae V,E_{n+k}\ad \ae(\bar\nabla_X J_k)W,Y\ad\\
& & \,\,\,\,\,+\frac{1}{2}\sum_k\ae Y,E_{n+k}\ad \ae(\bar\nabla_X
J_k )W,V\ad
+\frac{1}{2}\sum_k\ae W, E_{n+k}\ad \ae(\bar\nabla_Y J_k)V,X\ad\\
& &\,\,\,\,\, -\frac{1}{2}\sum_k\ae V,E_{n+k}\ad \ae(\bar\nabla_Y
J_k)W,X\ad-\frac{1}{2} \sum_k\ae X,E_{n+k}\ad \ae(\bar\nabla_Y J_k
)W,V\ad
\end{eqnarray*}
and
\begin{eqnarray*}
& & Q_2(X, Y, V, W)\\
& &\,\, =-\frac{1}{4}\sum_{k,l}\langle E_{n+k},W\rangle\langle
E_{n+l},V\rangle\langle J_k Y,J_l
X\rangle+\frac{1}{4}\sum_{k,l}\langle E_{n+k},W\rangle\langle
E_{n+l},X\rangle\langle J_k Y,J_l
V\rangle\\
& & \,\,\,\,-\frac{1}{4}\sum_{k,l}\langle E_{n+k},Y\rangle\langle
E_{n+l},V\rangle\langle J_k W,J_l
X\rangle+\frac{1}{4}\sum_{k,l}\langle E_{n+k},Y\rangle\langle
E_{n+l},X\rangle\langle J_k W,J_l
V\rangle\\
& & \,\,\,\,+\frac{1}{4}\sum_{k,l}\langle E_{n+k},W\rangle\langle
E_{n+l},V\rangle\langle J_k X,J_l
Y\rangle-\frac{1}{4}\sum_{k,l}\langle E_{n+k},W\rangle\langle
E_{n+l},Y\rangle\langle J_k X,J_l
V\rangle\\
& & \,\,\,\,+\frac{1}{4}\sum_{k,l}\langle E_{n+k},X\rangle\langle
E_{n+l},V\rangle\langle J_k W,J_l
Y\rangle-\frac{1}{4}\sum_{k,l}\langle E_{n+k},X\rangle\langle
E_{n+l},Y\rangle\langle J_k W,J_l V\rangle.
\end{eqnarray*}

An important relation between $\lambda$ and $Q$ is given by the
following lemma
\begin{lemma}\label{Q}
The components $Q^b_a$ of $Q$ are given by the $2$-forms
\begin{equation}
\label{q} Q^a_b := Q(\,\cdot \,, \, \cdot\,, e_b, e_a)=
\big(\dd\lambda+\lambda\wedge\omega+\omega\wedge\lambda-\lambda\wedge\lambda\big)\,^a_b.
\end{equation}
\end{lemma}

\noindent \emph{Proof.} Denoting the right hand side in  (\ref{q})
by $\Lambda^a_b$ and expanding it, it results that
\begin{eqnarray}
& & -2\Lambda^a_d = \sum_{c}\sum_k \big((\dd U^k_a-\sum_b
U^k_b\omega^b_a) u^k_{dc}-(\dd U^k_d-\sum_b U^k_b \omega^b_d)
u^k_{ac}-(\dd U^k_c-\sum_b U^k_b\omega^b_c)
u^k_{ad}\nonumber\\
& & \,\,\,\,\,\,\,\,+\,U^k_a (\dd u^k_{dc}-\sum_b u^k_{db}\omega^b_c-\sum_b u^k_{bc}\omega^b_d)\nonumber\\
& & \,\,\,\,\,\,\,\,-\,U^k_d (\dd u^k_{ac}-\sum_b u^k_{ab}\omega^b_c-\sum_b u^k_{bc}\omega^b_a)\nonumber\\
& & \,\,\,\,\,\,\,\,-\,U^k_c (\dd u^k_{ad}-\sum_b u^k_{ab}\omega^b_d-\sum_b u^k_{bd}\omega^b_a)\nonumber\\
& &\label{dlambda2} \,\,\,\,\,\,\,\,+\, U^k_a \sum_b
u^k_{bc}\lambda^b_d - \sum_b U^k_b \lambda^b_d u^k_{ac} -
U^k_c\sum_b u^k_{ab}\lambda^b_d\big)\wedge\omega^c.
\end{eqnarray}
The covariant derivative of the $(0,2)$ tensor $J_{k}$ has
components given in terms of the frame $\{e_a\}_{a=1}^{n+n'}$ by
\begin{eqnarray}
\label{Ju} \bar \nabla J_{k}(e_a, e_b)= \dd
u^k_{ab}-u^k_{db}\omega^d_a - u^k_{ad}\omega^d_b=:\bar\nabla
u^k_{ab}.
\end{eqnarray}
Using (\ref{u}) and (\ref{Ju}), one gets
\begin{eqnarray*}
& & -2\Lambda^a_d
=\,\sum_{k,c,c'}\big(-\frac{1}{2}u^k_{c'a}
u^k_{dc}+\frac{1}{2}u^k_{c'd} u^k_{ac}+\frac{1}{2}u^k_{c'c}
u^k_{ad}\big)\omega^{c'}\wedge
\omega^c\nonumber\\
& & \,\,\,\,\,\,\,\,+\,\sum_{k,c}\big(U^k_a \bar\nabla
u^k_{dc}-U^k_d \bar\nabla u^k_{ac}
 -U^k_c \bar\nabla u^k_{ad}\big)\wedge\omega^c\\
& &\,\,\,\,\,\,\,\,+\,\sum_{k,c}\big( U^k_a \sum_b
u^k_{bc}\lambda^b_d - \sum_b U^k_b \lambda^b_d u^k_{ac} -
U^k_c\sum_b u^k_{ab}\lambda^b_d\big)\wedge\omega^c.
\end{eqnarray*}
The last three terms may be calculated using that for $1\le k\le
n'$, $1\le a\le n+n'$, one has
\begin{equation}
\label{useful} \dd U^k_a -\sum_c U^k_c \omega^c_a + \sum_c U^k_c
\lambda^c_a=0.
\end{equation}
For proving (\ref{useful}), using (\ref{L}), one computes
\begin{eqnarray*} \sum_c U^k_c \lambda^c_a & =
& \sum_c U^k_c
L(\cdot,e_c,e_a) =  L(\cdot, E_{n+k},e_a)\\
&=& -\frac{1}{2}\sum_l \Big(\ae J_l e_a,\cdot\ad \ae
E_{n+k},E_{n+l}\ad-\ae J_l E_{n+k},\cdot\ad \ae
e_a,E_{n+l}\ad\\
& &\,\,+ \ae J_l E_{n+k},e_a\ad \ae \cdot, E_{n+l}\ad\Big)\\
& = & -\frac{1}{2} \big(\ae  J_k e_a,\cdot\ad -\sum_l \ae J_l
E_{n+k},\cdot\ad \ae e_a, E_{n+l}\ad -\sum_l \ae J_l E_{n+k},e_a\ad
\ae \cdot, E_{n+l}\ad\big).
\end{eqnarray*}
However, given any vector field $V$ in $N$, one has
\begin{eqnarray*}
\ae  J_l E_{n+k},V\ad = \sum_{r,s}\ae E_{n+k},E_r\ad \ae V,
E_s\ad\sigma_{rs}^{n+l}=\sum_s\ae V, E_s\ad \sigma_{n+k,s}^{n+l}=0.
\end{eqnarray*}
Therefore, one concludes that
\begin{eqnarray*}
\sum_c U^k_c \lambda^c_a  =   -\frac{1}{2}\ae J_k e_a,\cdot\ad =
-\frac{1}{2}\sum_b u^k_{ab}\omega^b =  \dd U^k_a - \sum_c
U^k_c\omega^c_a,
\end{eqnarray*}
as desired. This proves (\ref{useful}). Thus, we may write
\begin{eqnarray*}
& & -2\Lambda^a_d
=\sum_{k,c,c'}\big(-\frac{1}{2}u^k_{c'a} u^k_{dc}+\frac{1}{2}u^k_{c'
c} u^k_{ad}\big)\omega^{c'}\wedge
\omega^c\nonumber\\
& & \,\,\,\,\,\,\,\,+\,\sum_{k,c}\big(U^k_a \bar\nabla
u^k_{dc}-U^k_d \bar\nabla u^k_{ac}
 -U^k_c \bar\nabla u^k_{ad}\big)\wedge\omega^c\\
& &\,\,\,\,\,\,\,\,+\sum_{k,c}\big( U^k_a \sum_b u^k_{bc}\lambda^b_d
- U^k_c\sum_b u^k_{ab}\lambda^b_d\big)\wedge\omega^c.
\end{eqnarray*}
Nevertheless, in view of (\ref{L}), it follows that
\begin{eqnarray*}
& &  U^k_a \sum_b u^k_{bc}\lambda^b_d  + U^k_c\sum_b
u^k_{ba}\lambda^b_d = \sum_b \big(\langle E_{n+k},e_a\rangle\langle
J_k e_b,e_c\rangle +\langle E_{n+k},e_c\rangle\langle J_k
e_b,e_a\rangle\big) L(\cdot, e_b,e_d)\\
& & \,\, =\frac{1}{2}\sum_l \langle E_{n+k},e_a\rangle\langle
E_{n+l},e_d\rangle\langle J_k
e_c,J_l\cdot\rangle-\frac{1}{2}\sum_l\langle
E_{n+k},e_a\rangle\langle E_{n+l},\cdot\rangle\langle J_k e_c,J_l
e_d\rangle\\
& & \,\,\,\,\,\,\,+\frac{1}{2}\sum_l\langle
E_{n+k},e_c\rangle\langle E_{n+l},e_d\rangle\langle J_k
e_a,J_l\cdot\rangle -\frac{1}{2}\sum_l\langle
E_{n+k},e_c\rangle\langle E_{n+l},\cdot\rangle\langle J_k e_a,J_l
e_d\rangle.
\end{eqnarray*}
Therefore, one concludes that
\begin{eqnarray*}
& & \Lambda^a_d
=\sum_{k,c,c'}\big(\frac{1}{4}u^k_{c'a}
u^k_{dc}-\frac{1}{4}u^k_{c'c} u^k_{ad}\big)\omega^{c'}\wedge
\omega^c\nonumber\\
& & \,\,\,\,\,\,\,\,\,\,\,\,-\frac{1}{2}\sum_{k,c,c'}\big(U^k_a
\bar\nabla_{c'} u^k_{dc}-U^k_d \bar\nabla_{c'} u^k_{ac}
 -U^k_c \bar\nabla_{c'} u^k_{ad}\big)\omega^{c'}\wedge\omega^c\\
& &\,\,\,\,\,\,\,\,\,\,\,\,-\sum_{k,l,c,c'}\big( \frac{1}{4}\langle
E_{n+k},e_a\rangle\langle E_{n+l},e_d\rangle\langle J_k e_c,J_l
e_{c'}\rangle-\frac{1}{4}\langle E_{n+k},e_a\rangle\langle
E_{n+l},e_{c'}\rangle\langle J_k e_c,J_l
e_d\rangle\\
& & \,\,\,\,\,\,\,\,\,\,\,\,+\frac{1}{4}\langle
E_{n+k},e_c\rangle\langle E_{n+l},e_d\rangle\langle J_k e_a,J_l
e_{c'}\rangle-\frac{1}{4}\langle E_{n+k},e_c\rangle\langle
E_{n+l},e_{c'}\rangle\langle J_k e_a,J_l
e_d\rangle\big)\omega^{c'}\wedge\omega^c,
\end{eqnarray*}
what finishes the proof of the lemma. \hfill $\square$

\vspace{3mm}

\noindent This lemma has the following consequence, which
characterizes geometrically the tensor $Q$.

\begin{proposition}\label{curvQ} The tensor $Q$ satisfies the equation
\begin{equation}
\label{Qad} Q=A^{-1}\Theta A
\end{equation}
where $\Theta = (\Theta^k_l)_{k,l=1}^{n+n'}$ are the curvature forms
defined in (\ref{curv-amb}).
\end{proposition}

\noindent \emph{Proof.} One has
\begin{eqnarray}
\label{omegaQ}
&& d\omega +\omega\wedge\omega\nonumber \\
&& \,\,=d(A^{-1}dA)+A^{-1}dA\wedge A^{-1}dA +d\lambda +\lambda\wedge\lambda +\lambda\wedge A^{-1}dA+A^{-1}dA\wedge\lambda\nonumber\\
&& \,\,= -A^{-1}dA\wedge A^{-1} dA +A^{-1}dA\wedge A^{-1}dA
+d(A^{-1}\theta A) +A^{-1}\theta\wedge\theta A\nonumber\\
& & \,\,\,\,+A^{-1}\theta\wedge dA +A^{-1}dA A^{-1}\wedge\theta A\nonumber\\
& & \,\,=dA^{-1}\wedge \theta A +A^{-1}\dd\theta A-A^{-1}\theta\wedge dA +A^{-1}\theta\wedge\theta A\nonumber\\
& & \,\,\,\,+A^{-1}\theta\wedge dA -dA^{-1}\wedge\theta A\nonumber\\
&& \,\,=A^{-1}(d\theta+\theta\wedge\theta)A=A^{-1}\Theta A.
\end{eqnarray}
On the other hand we have
\begin{equation}
\label{zce0} d(\omega-\lambda)+(\omega-\lambda)\wedge
(\omega-\lambda)=-A^{-1}dA\wedge A^{-1}dA + A^{-1}dA\wedge
A^{-1}dA=0,
\end{equation}
what implies that
\begin{equation}
\label{omegalambda} d\omega +\omega\wedge\omega =
d\lambda-\lambda\wedge\lambda
+\omega\wedge\lambda+\lambda\wedge\omega.
\end{equation}
Hence (\ref{omegaQ}) and (\ref{omegalambda}) give the desired
result. \hfill $\square$

\section{Isometric immersions into two-step nilpotent Lie
groups}\label{imm-nec}

From now on, we consider a simply connected Riemannian manifold
$M^m$. We denote $m'=n+n'-m$.

From the calculations above, we infer the following necessary
conditions for the existence of isometric immersions in $N$ with
prescribed second fundamental form. In the statement, $\bar R$
denotes the curvature tensor in $N$.

\begin{proposition} Let $f:M\to N$ 
be an isometric immersion.
Then, the Gauss, Ricci and Codazzi equations are given by
\begin{equation}
\label{gcr-i} \bar R(f_*X,f_*Y,V,W)=Q(f_*X,f_*Y,V,W), \quad
X,Y\in\Gamma(TM),\,\, V,W\in\Gamma(f^*TN).
\end{equation}
Moreover, the following additional equations are satisfied
\begin{equation}
\label{eq-add-i} \bar\nabla_X E_{n+k}=-\frac{1}{2} J_k X, \quad
X\in\Gamma(TM)
\end{equation}
for $k=1,\ldots, n'$.
\end{proposition}

\noindent \emph{Proof.}  After identifying $M$ and the immersed
submanifold $f(M)\subset N$, we consider an orthonormal frame
$\{e_a\}_{a=1}^{m+m'}$ defined in an ambient open neighborhood of an
arbitrary point in $M$.  This frame may be chosen adapted to the
immersion, that is, in such a way that, along points in $M$, the
first $m$ fields in this frame are tangent to $M$ and the other $m'$
ones are local sections of the normal bundle $TM_f^\perp$.

Let $A$ be given as above by (\ref{eE}). Then, 
the connection forms $\{\omega^a_b\}_{a,b=1}^{m+m'}$ satisfy
\begin{equation}
\label{gcr-ad} d\omega^a_b+\sum_c\omega^a_c\wedge\omega^c_b =
(A^{-1}\Theta A)^a_b,
\end{equation}
where $\Theta$ is given by (\ref{curv-amb}). Since the right-hand
side in (\ref{gcr-ad}) corresponds to the ambient curvature tensor
expressed in terms of the adapted frame, this equation corresponds
to Gauss, Codazzi and Ricci equations, respectively, as we may
easily verify considering suitable ranges of indices $a,b$. Hence,
(\ref{Qad}) in Proposition \ref{curvQ} implies (\ref{gcr-i}).

The equation (\ref{eq-add-i}) follows immediately from the preceding
discussion.
\hfill $\square$

\section{Existence of an adapted frame}\label{ex-frame}


We now consider a real Riemannian vector bundle $\mathcal{E}$ over
$M$ with rank $m'$ and the Whitney sum bundle
$\mathcal{S}=TM\oplus\mathcal{E}$. The metric in $\mathcal{S}$ is
also represented by $\langle\cdot,\cdot\rangle$. Let $\hat\nabla$
and $\hat R$ be respectively the compatible connection and curvature
tensor in $\mathcal{S}$.

We suppose that $\mathcal{S}$ is a trivial vector bundle and then we
fix a globally defined orthonormal frame $\hat E_1, \ldots, \hat
E_{n+n'}$ in $\mathcal{S}$. Hence, for any $k=1,\ldots, n'$, one
defines
\begin{equation}
\label{Jk} \ae \hat J_k V,W\ad =\sum_{l,r=1}^n \ae V,\hat E_l\ad \ae
W,\hat E_r\ad \sigma^{n+k}_{lr},\quad V,W\in \Gamma(\mathcal{S}),
\end{equation}
where the constants $\sigma^{n+k}_{lr}$ are given by (\ref{sigma0}).
It is obvious from the definition that
\begin{equation}
\label{kk0} \ae\hat J_k V,\hat E_{n+l}\ad = 0
\end{equation}
since $\sigma^{n+k}_{r,n+l}=0$.

Now, we define in terms of $\hat J_k$ tensors $\hat L$ and $\hat Q$
in $\mathcal{S}$ by
\begin{eqnarray}
\label{LA}\hat L(X,Y,V)& = &  -\frac{1}{2}\sum_{k=1}^{n'}\ae \hat
J_k V,X\ad \ae Y,\hat E_{n+k}\ad+\frac{1}{2}\sum_{k=1}^{n'}\ae \hat
J_k Y,X\ad \ae V,\hat E_{n+k}\ad
\nonumber\\
& &\,-\frac{1}{2}\sum_{k=1}^{n'}\ae \hat J_{k}Y,V\ad \ae X,\hat
E_{n+k}\ad, \quad X,Y,V\in \Gamma(\mathcal{S})
\end{eqnarray}
and for $X,Y\in\Gamma(TM)$ and $V,W\in \Gamma(\mathcal{S})$,
\begin{eqnarray}
\label{QA} \hat Q(X,Y,V,W)=\hat Q_1(X,Y,V,W)+\hat Q_2(X,Y,V,W),
\end{eqnarray}
where
\begin{eqnarray*}
\label{QA1}
&  & \hat Q_1(X, Y, V, W)\nonumber\\
& &\,\,\,\,\,=\frac{1}{4}\ae \hat J_k X,W\ad\ae \hat J_k V,Y\ad
+\frac{1}{2}\ae \hat J_k Y,X\ad \ae \hat J_k W,V\ad-\frac{1}{4}\ae
\hat J_k Y,W\ad \ae
\hat J_k V,X\ad\\
& & \,\,\,\,\,-\frac{1}{2}\sum_k\ae W,\hat E_{n+k}\ad
\ae(\hat\nabla_X \hat J_k)V,Y\ad
+\frac{1}{2}\sum_k\ae V,\hat E_{n+k}\ad \ae(\hat\nabla_X \hat J_k)W,Y\ad\\
& & \,\,\,\,\,+\frac{1}{2}\sum_k\ae Y,\hat E_{n+k}\ad
\ae(\hat\nabla_X \hat J_k )W,V\ad
+\frac{1}{2}\sum_k\ae W, \hat E_{n+k}\ad \ae(\hat\nabla_Y \hat J_k)V,X\ad\\
& &\,\,\,\,\, -\frac{1}{2}\sum_k\ae V,\hat E_{n+k}\ad
\ae(\hat\nabla_Y \hat J_k)W,X\ad-\frac{1}{2} \sum_k\ae X,\hat
E_{n+k}\ad \ae(\hat\nabla_Y \hat J_k )W,V\ad
\end{eqnarray*}
and
\begin{eqnarray*}
\label{QA2}
&  & \hat Q_2(X, Y, V, W)\nonumber\\
& &\,\,\,\,\, = -\frac{1}{4}\sum_{k,l}\langle \hat
E_{n+k},W\rangle\langle \hat E_{n+l},V\rangle\langle \hat J_k Y,\hat
J_l X\rangle+\frac{1}{4}\sum_{k,l}\langle\hat
E_{n+k},W\rangle\langle\hat E_{n+l},X\rangle\langle \hat J_k Y,\hat
J_l
V\rangle\\
& & \,\,\,\,-\frac{1}{4}\sum_{k,l}\langle \hat
E_{n+k},Y\rangle\langle\hat E_{n+l},V\rangle\langle \hat J_k W,\hat
J_l X\rangle+\frac{1}{4}\sum_{k,l}\langle\hat
E_{n+k},Y\rangle\langle \hat E_{n+l},X\rangle\langle \hat J_k W,\hat
J_l
V\rangle\\
& & \,\,\,\,+\frac{1}{4}\sum_{k,l}\langle \hat
E_{n+k},W\rangle\langle \hat E_{n+l},V\rangle\langle\hat J_k X,\hat
J_l Y\rangle-\frac{1}{4}\sum_{k,l}\langle \hat
E_{n+k},W\rangle\langle \hat E_{n+l},Y\rangle\langle \hat J_k X,\hat
J_l
V\rangle\\
& & \,\,\,\,+\frac{1}{4}\sum_{k,l}\langle \hat
E_{n+k},X\rangle\langle \hat E_{n+l},V\rangle\langle \hat J_k W,\hat
J_l Y\rangle-\frac{1}{4}\sum_{k,l}\langle \hat
E_{n+k},X\rangle\langle\hat E_{n+l},Y\rangle\langle \hat J_k W,\hat
J_l V\rangle.
\end{eqnarray*}
We then suppose that
\begin{equation}
\label{GCR} \ae\hat R(X,Y)V,W\ad= \hat Q(X,Y,V,W), \quad X,Y \in
\Gamma(TM), \,\, V,W \in \Gamma(\mathcal{S}).
\end{equation}
We also assume the following condition
\begin{equation}
\label{additional-0} \hat \nabla_X \hat E_{n+k} = -\frac{1}{2}\hat
J_k X,\quad X\in \Gamma(TM),\quad k=1,\ldots, n'.
\end{equation}
The connection in $\mathcal{S}$ induces connections $\nabla$ in $M$
and $\nabla^{\mathcal{E}}$ in $\mathcal{E}$. More precisely,
defining $II\in \Gamma(T^*M\otimes T^*M\otimes\mathcal{E})$ by
\begin{equation}
\hat\nabla_X Y = \nabla_X Y +II(X,Y), \quad X,Y \in \Gamma(TM)
\end{equation}
and defining, for $V\in \Gamma(\mathcal{E})$,
\begin{equation}
\ae S_V(X),Y\ad =\ae II(X,Y),V\ad,
\end{equation}
one obtains
\begin{equation}
\hat\nabla_X V = -S_V X + \nabla^{\mathcal{E}}_X V.
\end{equation}
In terms of the decomposition $\hat E_{n+k}= T_k + N_k$, $T_k\in
\Gamma(TM)$, $N_k\in \Gamma(\mathcal{E})$, the condition
(\ref{additional-0}) becomes
\begin{eqnarray}
\label{additional} \nabla_X T_k -S_k(X) + \nabla^{\mathcal{E}}_X N_k
+ II(T_k,X) = -\frac{1}{2}\hat J_k (X), \quad X\in \Gamma(TM),
\end{eqnarray}
where $S_k = S_{N_k}$.

\begin{definition}
Given a connected simply connected open subset $M'\subset M$, we fix
a map $\hat U\in C^\infty(M',\mathbb{R}^{n'(n+n')})$. A frame $e:M'
\to \mathcal{F}(\mathcal{S})$ with components
\[
e_1,\ldots e_m, e_{m+1},\ldots e_{m+m'}
\]
is \emph{admissible} if the first $m$ sections are vector fields in
$M'$ and the last $m'$ ones are sections in $\mathcal{E}$ and,
moreover,  if it holds  that
\begin{equation}
\ae\hat  E_{n+k},e_a\ad = \hat U^k_a, \quad 1\le k\le n'.
\end{equation}
In particular, this implies that
\begin{equation}
\label{Tk} \ae T_k, e_i\ad = \ae \hat E_{n+k},e_i\ad = \hat U^k_i
\end{equation}
for $i=1,\ldots, m$ and
\begin{equation}
\label{Nk} \ae N_k, e_\alpha\ad = \ae \hat E_{n+k},e_\alpha\ad =
\hat U^k_\alpha
\end{equation}
for $\alpha=m+1,\ldots, m+m'$. The transition map from the frame
$\{\hat E_k\}_{k=1}^{n+n'}$ to an admissible frame
$\{e_a\}_{a=1}^{m+m'}$ is given by an \emph{admissible map}, that
is, if
\begin{equation}
\label{ad-frame} e_a = \sum_{k=1}^{n+n'}\hat E_k A^k_a,
\end{equation}
then $A$ is if the form
\begin{eqnarray}
\label{admissible}
A(x)= \left(\begin{array}{cc} * \\
\hat U(x)
\end{array}\right),
\end{eqnarray}
where the block $\hat U(x)$ corresponds to the last $n'$ lines.
\end{definition}

\noindent We denote by
\begin{equation}
\omega^1,\ldots, \omega^m, \omega^{m+1},\ldots, \omega^{m+m'}
\end{equation}
the real-valued $1$-forms   dual to the frame
$\{e_a\}_{a=1}^{m+m'}$. The Riemannian connection $\hat\nabla$ is
given  in terms of this  frame by the matrix
$\omega=(\omega^a_b)_{a,b=1}^{n+n'}$. Hence,  the first structural
equation is written as
\begin{equation}
\label{first-formal} \dd \omega^a + \sum_b \omega^a_b\wedge \omega^b
= 0, \quad \omega^a_b = -\omega^b_a.
\end{equation}
Regarding each $\hat J_k$ as  $(0,2)$ tensor, we write them locally
as
\begin{equation}
\hat J_k= \sum_{a,b} \hat u^k_{ab}\omega^a \otimes \omega^b, \quad k
=1,\ldots, n'.
\end{equation}
Thus in local terms the equation (\ref{additional-0}) is rewritten
as
\begin{equation}
\label{add-local} \sum_{k}\big( \dd \hat U^k_a-\sum_b \hat
U^k_b\omega^b_a\big)=\frac{1}{2}\sum_k \hat u^k_{ab}\,\omega^b.
\end{equation}
The local expression for $\hat L$ is given by the $1$-forms
\begin{equation}
\hat\lambda^a_b =\hat L(\,\cdot\, , e_a, e_b).
\end{equation}
Following the calculations in Section \ref{ctensor}, we conclude
that
\begin{equation}
\hat\lambda^a_b =-\frac{1}{2}\sum_{c=1}^{n+n'}\sum_{k=1}^{n'}
\big(\hat U^k_a \hat u^k_{bc}-\hat U^k_b\hat u^k_{ac}+\hat U^k_c
\hat u^k_{ba}\big)\omega^c.
\end{equation}
The local expression for $\hat Q$ is given by the $2$-forms
\begin{equation}
\hat Q_d^a = \hat Q(\cdot, \cdot, e_d,e_a).
\end{equation}
One notices that the hypothesis (\ref{GCR}) is rephrased in terms of
these forms as
\begin{equation}
\label{eqtn-for-omega} \dd\omega^a_b+ \sum_c\omega^a_c\wedge
\omega^c_b = \hat Q^a_b.
\end{equation}
We may verify proceeding as in the  proof of the Lemma \ref{Q} and
using (\ref{kk0}) that
\begin{equation}
\label{eqtn-for-lambda} \hat Q^a_d := \hat Q(\,\cdot, \,, \,
\cdot\,, e_d, e_a)=
\big(\dd\hat\lambda+\hat\lambda\wedge\omega+\omega\wedge\hat\lambda-\hat\lambda\wedge\hat\lambda\big)\,^a_d.
\end{equation}
Combining equations (\ref{eqtn-for-omega}) and
(\ref{eqtn-for-lambda}), one deduces that $\hat\omega:=\omega-\hat
\lambda$ satisfies  the zero curvature equation
\begin{equation}
\label{zce} \dd\hat\omega + \hat\omega \wedge \hat\omega =0.
\end{equation}
A suitable version of (\ref{useful}) allows us to claim that
\begin{equation}
\label{useful2} \dd \hat U^k_a -\sum_b \hat U^k_b\omega^b_a+ \sum_b
\hat U^k_b\hat \lambda^b_a=0.
\end{equation}
We then prove the following result.

\begin{proposition}
\label{frame} Assume that (\ref{GCR}) and (\ref{additional-0}) hold.
Let $M'\subset M$ be a connected simply connected open subset. Then
there exists an admissible map $A\in
C^\infty(M',\textsc{SO}_{n+n'})$ so that
\begin{equation}
\label{gauge} A^{-1}\dd A = \omega-\hat\lambda
\end{equation}
with initial condition $A(x_0)=\textrm{Id}$, for a given $x_0\in
M'$.
\end{proposition}

\noindent \emph{Proof.}  We want to assure the existence of an
admissible map so that
\begin{equation}
\label{eqtn-for-A} A^{-1}\dd A = \hat\omega,
\end{equation}
where $\hat\omega=\omega-\hat\lambda$.

If we denote by $\mu:\textsc{M}_{n+n'}\mathbb{R}\to
\mathbb{R}^{n'(n+n')}$ the projection on the last $n'$ lines, the
condition (\ref{ad-frame}) means that $\mu(A(x))=\hat U(x)$. The set
of admissible maps define a submanifold of $M\times
\textsc{SO}_{n+n'}$, namely
\begin{equation}
\mathcal{U} =\bigg\{(x,A):A=\left(\begin{array}{cc} * \\
\hat U(x)\end{array} \right)\bigg\},
\end{equation}
whose  tangent space  at a point $(x,A)$ is
\begin{equation}
T_{(x,A)}\mathcal{U}=\bigg\{(v,\textrm{\bf B}):\textrm{\bf B}=
\left(\begin{array}{cc} * \\
\dd \hat U(x)\cdot v\end{array} \right)\bigg\}.
\end{equation}

Let $\bar\omega\in
\Lambda^1(\textsc{SO}_{n+n'},\mathfrak{so}_{n+n'})$ be the
Maurer-Cartan form in $\textsc{SO}_{n+n'}$. Thus, the equation
(\ref{eqtn-for-A}) is written as
\begin{equation}
\label{Dar-der} \hat\omega=A^*\bar\omega.
\end{equation}
For solving this equation, we define a $1$-form $\Upsilon$ in
$M'\times \textsc{SO}_{n+n'}$ with values on $\mathfrak{so}_{n+n'}$
by
\begin{equation}
\Upsilon = \pi^*_1\hat\omega -\pi^*_2\bar\omega,
\end{equation}
where $\pi_1: M\times \textsc{SO}_{n+n'}\to M$ and $\pi_2:M\times
\textsc{SO}_{n+n'}\to \textsc{SO}_{n+n'}$ are the natural
projections. We then define the distribution
$\mathcal{D}=\ker\Upsilon$ on $\mathcal{U}$. More precisely
\begin{equation}
\label{DD} (v,{\bf B}) \in \mathcal{D}_{(x,A)} \quad \textrm{ if and
only if}\quad \hat\omega_x(v)=\bar\omega_A({\bf B})
\end{equation}
In order to prove that (\ref{DD}) defines a distribution we must
verify that $\ker \Upsilon$ has constant rank. We begin by proving
that the differential of $\pi_1$ restricted to $\mathcal{D}_{(x,A)}$
is a monomorphism. In fact, if $\pi_{1_*} (v,{\bf B}) = 0$ for some
$(v, {\bf B})\in \mathcal{D}_{(x,A)}$ then $v=0$. Since $0= \hat
\omega_x (v) = \bar\omega_A({\bf B})$, it follows that ${\bf B}=0$.
Therefore,
\[
\dim\ker\Upsilon_{(x,\A)}\le m.
\]
Now, given $(v,{\bf B})\in T_{(x,A)}\mathcal{U}$ we have
\begin{eqnarray*}
& & \mu\big(A\Upsilon_{(x,A)}(v,{\bf B})\big) =
\mu\big(A\hat\omega_x(v)-A\bar\omega_{A}({\bf
B})\big)=\mu\big(A\hat\omega_x(v)-AA^{-1}\cdot {\bf
B}\big)\\
& & \,\,=\mu\big(A\big)\hat\omega_x(v)-\mu\big({\bf B}\big)=\hat
U\hat\omega_x(v)-\dd \hat U_x(v)
=0
\end{eqnarray*}
where  in the last equality we used equation (\ref{useful2}). We
then had verified that
\[
\textrm{Im}\Upsilon_{(x,A)}\subset
\{\mathcal{B}\in\mathfrak{so}_{n+n'}: \mu(A\mathcal{B})=0\}
\]
Thus, if $\mathcal{B}\in\textrm{Im}\Upsilon_{(x,A)}$ then
$\mathcal{B}=\bar\omega_{A}({\bf B})$ for some $\bf B$ tangent to
$A$ such that $\mu({\bf B})=0$. This means that
\[
\textrm{Im}\Upsilon_{(x,A)}\subset
\bar\omega_{A}\big(\ker\mu_{A}\big)
\]
where $\ker \mu_{A}=\{{\bf B}\in T_{A}\textsc{SO}_{n+n'}: \mu({\bf
B})=0\}$. Since $\bar\omega_A$ is an isomorphism, it follows that
$\bar\omega_{A}\big(\ker\mu_{A}\big)$ and $\ker\mu_A$ have same
dimension. Thus,
\[
\dim \ker \Upsilon_{(x,A)}\ge m.
\]
Hence, $\mathcal{D}_{(x,A)}$ is $m$ dimensional, for all $x\in M',
A\in \mathcal{U}$.

Now we verify the integrability of $\mathcal{D}$. The zero curvature
equation (\ref{zce}) implies that
\begin{eqnarray*}
\dd \Upsilon &=& \dd \hat \omega - \dd \bar\omega = \hat \omega
\wedge \hat \omega - \bar\omega \wedge \bar\omega =(\bar\omega +
\Upsilon)\wedge (\bar\omega + \Upsilon)- \bar\omega \wedge \bar\omega\\
& = &\bar\omega \wedge \Upsilon+ \Upsilon \wedge \bar\omega.
\end{eqnarray*}
Thus if one calculates $\dd \Upsilon$ at some vector $(v,{\bf B})\in
\mathcal{D}_{(x,A)}$ one obtains $\Upsilon (v,{\bf B}) =0$ and then
$\dd \Upsilon (v,{\bf B})=0$ too. So the ideal $\ker \Upsilon$ is
differential and then the distribution $\mathcal{D}$ is integrable.

Since $\pi$ is a local diffeomorphism  between the simply connected
domain $M'$ and the integral leaf of $\mathcal{D}$ passing through
$(x_0, \textrm{Id})$, a standard monodromy reasoning implies that
 this leaf as the graph $x\mapsto A(x)$ of a certain
map $A\in C^\infty(M',\textsc{SO}_{n+n'})$ which by definition
satisfies (\ref{ad-frame}) and (\ref{Dar-der}). \hfill $\square$

\vspace{3mm}

\noindent Given an admissible map $A:M'\to \textsc{SO}_{n+n'}$
solving (\ref{gauge}), one defines a frame $\{e_a\}_{a=1}^{m+m'}$ in
$\mathcal{S}$ along $M'$ by (\ref{ad-frame}).  The corresponding
sets of dual $1$-forms are related by
\begin{equation}
\label{theta-A} \hat\theta^k = \sum_{a=1}^{m+m'}A^k_a \omega^a.
\end{equation}
It stems from (\ref{Jk}) that the local expression for $\hat J_k$ in
the frame $\{e_a\}_{a=1}^{n+n'}$ is
\[
-\frac{1}{2}\hat u^k_{ab}=-\frac{1}{2}\ae \hat J_k e_a, e_b\ad =
\sum_{l,r}\ae e_a, \hat E_l\ad \ae e_b, \hat E_r\ad
\sigma^{n+k}_{lr}=\sum_{l,r} A^l_a A^r_b \sigma^{n+k}_{lr}.
\]
We then define
\begin{equation}
\label{Theta-A} \hat\theta^k_l
=\frac{1}{2}\sum_r\tau^k_{lr}\hat\theta^r =
\frac{1}{2}\sum_a\sum_r\tau^k_{lr} A^r_a \omega^a.
\end{equation}
In view of these facts, we are able to restate Proposition \ref{l}
in the current context.

\begin{proposition}
\label{l-A} The admissible frame obtained above as solution of the
equation (\ref{gauge}) satisfies
\begin{equation}
\label{l-theta} \hat\lambda=A^{-1}\hat\theta A,
\end{equation}
where $\hat\theta = (\hat\theta^k_l)_{k,l=1}^{n+n'}$ is defined in
(\ref{Theta-A}).
\end{proposition}

 \noindent{\it Proof.} It suffices to mimic the proof of Proposition \ref{l} in Section \ref{ctensor}.
\hfill $\square$

\vspace{3mm}

\noindent We finally define the following $2$-forms
\begin{equation}
\label{TTheta-A} \hat\Theta^k_l =
\frac{1}{4}\sum_{s,t}\big(\tau^k_{lr}\tau^r_{st}+\tau^k_{rs}\tau^r_{lt}\big)\hat\theta^s\wedge\hat\theta^t
=\frac{1}{4}\sum_{a,b}\sum_{s,t}\big(\tau^k_{lr}\tau^r_{st}+\tau^k_{rs}\tau^r_{lt}\big)A^s_a
A^t_b \omega^a\wedge\omega^b.
\end{equation}
Then we are able to prove the following result.

\begin{proposition}
\label{Q-A} The admissible frame defined above as solution of the
equation (\ref{gauge}) satisfies
\begin{equation}
\label{Q-theta} \hat Q = A^{-1}\hat\Theta A,
\end{equation}
where $\hat\Theta = (\hat\Theta^k_l)_{k,l=1}^{n+n'}$ is defined in
(\ref{TTheta-A}).
\end{proposition}

\noindent {\it Proof. } From (\ref{gauge}) and (\ref{l-theta}) it
follows that
\begin{eqnarray*} \dd \hat\lambda & = &  \dd A^{-1}\wedge\hat\theta
A + A^{-1}\dd \hat\theta A -
A^{-1}\hat\theta \wedge \dd A\\
&=& -A^{-1}\dd A \wedge A^{-1}\hat\theta A +A^{-1}\dd\hat\theta A
-A^{-1}\hat\theta A \wedge A^{-1}\dd A\\
& = & -(\omega-\hat\lambda)\wedge \hat\lambda +A^{-1}\dd \hat\theta
A -
\hat\lambda\wedge (\omega-\hat\lambda)\\
&=& 2 \hat\lambda\wedge\hat\lambda-\omega \wedge \hat\lambda
-\hat\lambda\wedge \omega+A^{-1}\dd \hat\theta A.
\end{eqnarray*}
Therefore, in view of (\ref{eqtn-for-lambda}), we conclude that
\begin{eqnarray}
\hat Q &=&
\dd\hat\lambda-\hat\lambda\wedge\hat\lambda+\omega\wedge\hat\lambda
+\hat\lambda\wedge \omega = \hat\lambda\wedge \hat\lambda +A^{-1}\dd
\hat\theta A\nonumber
\\
&=&
A^{-1}\hat\theta A\wedge A^{-1}\hat\theta A + A^{-1}\dd\hat\theta A\nonumber \\
\label{Qt}&  = & A^{-1}\big(\dd\hat\theta +\hat\theta \wedge
\hat\theta\big)A.
\end{eqnarray}
However, it follows from (\ref{Theta-A}), (\ref{first-formal}) and
(\ref{gauge}) that
\begin{eqnarray*}
\dd\hat\theta^k_l &=& \frac{1}{2}\sum_a\sum_r\tau^k_{lr}(\dd A^r_a
\wedge \omega^a + A^r_a \dd
\omega^a)=\frac{1}{2}\sum_{a,b}\sum_r\tau^k_{lr}(\dd A^r_b \wedge
\omega^b
- A^r_a \omega^a_b \wedge \omega^b)\\
&=&\frac{1}{2}\sum_{a,b}\sum_r\tau^k_{lr}(\dd A^r_b-A^r_a\omega^a_b)
\wedge
\omega^b \\
&  =& -\frac{1}{2}\sum_b\sum_r\tau^k_{lr} (A\hat\lambda)^r_b \wedge
\omega^b.
\end{eqnarray*}
However $A\hat\lambda = AA^{-1}\hat\theta A = \hat\theta A$. Hence,
one gets
\begin{eqnarray*}
\dd\hat\theta^k_l &=& -\frac{1}{2}\sum_b\sum_r\tau^k_{lr}(\hat\theta
A)^r_b\wedge \omega^b
=-\frac{1}{2}\sum_b\sum_{r,s}\tau^k_{lr}\hat\theta^r_s \wedge
A^s_b\omega^b= -\frac{1}{2}\sum_{r,s}\tau^k_{lr}\hat\theta^r_s
\wedge
\hat\theta^s\\
&=&-\frac{1}{4}\sum_{r,s,t}\tau^k_{lr}\tau^r_{st}\hat\theta^t \wedge
\hat\theta^s.
\end{eqnarray*}
On the other hand, one has
\begin{eqnarray*}
& & \sum_r\hat\theta^k_r \wedge \hat\theta^r_l =
\frac{1}{4}\sum_{r,s,t}\tau^k_{rs}\tau^r_{lt}\hat\theta^s\wedge\hat\theta^t.
\end{eqnarray*}
Therefore, one concludes that
\begin{equation}
\label{2nd-Formal} \dd\hat\theta + \hat\theta\wedge \hat\theta =
\hat\Theta.
\end{equation}
Gathering (\ref{Qt}) and (\ref{2nd-Formal}) we finish the proof.
\hfill $\square$

\section{Proof of the Theorem}\label{ex-imm}

\emph{Part a.} In view of the hypothesis in Theorem \ref{main},
Proposition \ref{frame} implies that there exists an admissible map
$A:M\to \textsc{SO}_{n+n'}$ which solves (\ref{gauge}) and satisfies
(\ref{l-theta}) and (\ref{Q-theta}) for
$\{\hat\theta^k\}_{k=1}^{n+n'}$ and
$\{\hat\theta^k_l\}_{k,l=1}^{n+n'}$ defined in (\ref{Theta-A}) and
(\ref{TTheta-A}), respectively.

We fix in $\mathfrak{n}=\mathbb{R}^{n+n'}$ the orthonormal frame
$\{\bar e_k=\omega_{\mathfrak{n}}(E_k)\}_{k=1}^{n+n'}$. We then
define the following $1$-form on $M\times N$ with values on
$\mathfrak{n}$
\[
\Pi =\pi_{N}^*\,\omega_{\mathfrak{n}}
-\sum_{k=1}^{n+n'}\sum_{a=1}^{m+m'} \bar e_k
(A^k_a\circ\pi_M)\pi_M^*\omega^a,
\]
where $\pi_{N}:M\times N\to N$ and $\pi_{M}:M\times N\to M$ are the
canonical projections. 
We then  consider the distribution $\mathcal{P} = \ker \Pi$ on
$M\times N$. Thus, using (\ref{first-formal}) and (\ref{gauge}), we
calculate (omitting projections)
\begin{eqnarray*}
 \dd \Pi &=&  \dd \omega_{\mathfrak{n}}-\sum_{a,k} \bar e_k \, \dd A^k_a \wedge
\omega^a -\sum_{a,k}\bar e_k\,
A^k_a\,\dd \omega^a\\
&  = &-\frac{1}{2} [\omega_{\mathfrak{n}}, \omega_{\mathfrak{n}}]
-\sum_{a,k}\bar e_k \,(A\hat \omega)^k_a\wedge \omega^a
+\sum_{a,c,k}\bar e_k A^k_a \,\omega^a_c \wedge
\omega^c\\
& = & -\frac{1}{2} [\Pi +\sum_k \bar e_k \hat\theta^k, \Pi +
\sum_l\bar e_l \hat\theta^l]-\sum_{a,k}\bar e_k \,(A\hat
\omega)^k_a\wedge \omega^a +\sum_{a,c,k}\bar e_k A^k_a \,\omega^a_c
\wedge \omega^c.
\end{eqnarray*}
Hence, one has
\begin{eqnarray*}
\dd\Pi &  = & -\frac{1}{2}[\Pi, \Pi]-\frac{1}{2}[\Pi, \sum_k\bar e_k
\hat\theta^k]-\frac{1}{2}[\sum_l\bar e_l\hat\theta^l,\Pi] -
\frac{1}{2}\sum_{k,l}[\bar e_k \hat\theta^k
, \bar e_l \hat\theta^l]\\
& &\,\, -\sum_{a,k}\bar e_k \, (A\omega)^k_a \wedge \omega^a
+\sum_{a,k}\bar e_k \,(A\hat\lambda)^k_a \wedge
\omega^a+\sum_{a,c,k}\bar e_k A^k_a \,\omega^a_c \wedge \omega^c.
\end{eqnarray*}
Thus  considering equality modulo $\Pi$ it follows that
\begin{eqnarray*}
\dd \Pi &=&
-\frac{1}{2}\sum_{k,l}\hat\theta^k\wedge \hat\theta^l \,[\bar
e_k,\bar e_l] -\sum_{a,c,k}\bar e_k\, A^k_c\omega^c_a\wedge \omega^a
+ \sum_{a,c,k}\bar e_k\, A^k_c\hat\lambda^c_a \wedge
\omega^a +\sum_{a,c,k}\bar e_k \,A^k_c\omega^c_a\wedge \omega^a\\
&  = & -\frac{1}{2}\sum_{k,l}\hat\theta^k\wedge \hat\theta^l \,[\bar
e_k,\bar e_l] + \sum_{a,c,k}\bar e_k A^k_c\hat\lambda^c_a \wedge
\omega^a.
\end{eqnarray*}
However using (\ref{l-theta}) one obtains
\begin{eqnarray*}
\dd \Pi &=& -\frac{1}{2}\sum_{k,l,r}\bar
e_r\sigma^r_{kl}\,\hat\theta^k\wedge \hat\theta^l +
\sum_{a,b,c}\sum_{k,l}\bar e_k\, A^k_c\hat\lambda^c_b (A^{-1})^b_l A^l_a\wedge \omega^a \\
& = & -\frac{1}{2}\sum_{k,l}\sum_r\bar
e_k\sigma^k_{rl}\,\hat\theta^r\wedge \hat\theta^l +\sum_{k,l}\bar
e_k \hat\theta^k_l\wedge \hat\theta^l= \sum_{k,l}\bar
e_k\big(\hat\theta^k_l
-\frac{1}{2}\sum_r\sigma^k_{rl}\hat\theta^r\big)\wedge
\hat\theta^l.
\end{eqnarray*}
Therefore  $\mathcal{P}$ is involutive since by (\ref{Theta-A}) one
has
\begin{eqnarray}
\label{connection_Lie_3}  \hat\theta^k_l =
\frac{1}{2}\sum_r\sigma^k_{rl}\hat\theta^r+\frac{1}{2}\sum_r\mu^k_{lr}\hat\theta^r,
\end{eqnarray}
where $\mu^k_{lr}=\sigma^l_{kr}+\sigma^r_{kl}$ satisfies
$\mu^k_{lr}=\mu^k_{rl}.$  This symmetry implies that
\[
\sum_{l}\big(\hat\theta^k_l -
\frac{1}{2}\sum_r\sigma^k_{rl}\hat\theta^r\big)\wedge \theta^l
=\frac{1}{2}\sum_{l,r} \mu^k_{lr}\hat\theta^r\wedge\hat\theta^l=0,
\]
what gives the integrability condition
\begin{equation*}
d\Pi = 0 \,\mod \Pi.
\end{equation*}
We may verify that an integral leaf through the identity $y_0$ in
$N$ is a graph over $M$. The function that graphics this leaf is an
isometric immersion $f:M\to N$ with initial condition, say,
$f(x_0)=y_0$, for a given point $x_0\in M$.

Indeed, given a tangent vector $(v,w)\in \mathcal{P}_{(x,y)}$  with
$y=f(x)$, we have $f_{*}(x)\cdot v=w$ and
\[
\omega_{\mathfrak{n}}(w) -\sum_{k=1}^{n+n'}\sum_{a=1}^{m+m'} \bar
e_k A^k_a(x)\omega^a(v)=0
\]
what yields after left translating both sides by $y$
\[
f_{*}(x)\cdot v =w= \sum_{k=1}^{n+n'}\sum_{a=1}^{m+m'}
E_k(f(x))A_a^k(x)\omega^a(v).
\]
Since $A(x)$ is an orthogonal matrix, we conclude that $f$ is an
isometric immersion and that
\[
e_a|_{f(x)}=\sum_{k=1}^{n+n'} E_k|_{f(x)}A^k_a(x),\quad 1\le a\le
m+m',
\]
defines an adapted frame along $f$ with corresponding dual co-frame
$\{\omega^a\}_{a=1}^{m+m'}$. Thus, it follows from
(\ref{first-formal}) that $\{\omega^a_b\}_{a,b=1}^{m+m'}$ are the
connection forms. Thus, (\ref{gauge}) and (\ref{l-theta}) imply that
$\{\hat\theta^k_l\}_{k,l=1}^{n+n'}$ are the connection forms in $N$
along $f$ with respect to the left-invariant frame
$\{E_k\}_{k=1}^{n+n'}$.  The equation (\ref{2nd-Formal}) assures
that $\{\hat\Theta^k_l\}_{k,l=1}^{n+n'}$ are the corresponding
curvature forms along $f$. Finally, (\ref{TTheta-A}) guarantees that
$\hat Q$ is the curvature form in $N$ at points of $f(M)$ associated
to the adapted frame $\{e_a\}_{a=1}^{m+m'}$.

 The choice of the initial
condition $f(x_0)=y_0$ is not a serious restriction, since an
isometric immersion with initial condition $y\in N$ is obtained
merely composing $f$ and the left translation by $yy_0^{-1}$.

\vspace{3mm}

\noindent \emph{Part b.} From (\ref{eq-add}) and (\ref{eq-add-2}) it
follows that that there exist local orthonormal frames
$\{e_a\}_{a=1}^{m+m'}$ and $\{\tilde e_a\}_{a=1}^{m+m'}$
respectively adapted to $f$ and $\tilde f$ such that the orthogonal
matrices
\begin{equation}
A^k_a =\langle e_a, E_k\rangle, \quad \tilde A^k_a = \langle \tilde
e_a, E_k\rangle
\end{equation}
satisfy
\begin{equation}
\label{mui} \mu(A)=\mu(\tilde A).
\end{equation}
Moreover, (\ref{IIi}) and (\ref{perpi}) imply that the connection
forms $\omega$ and $\tilde \omega$ for adapted frames along $f$ and
$\tilde f$ satisfy at corresponding points $\omega = \tilde\omega$.

Finally, (\ref{mui}), (\ref{u3}) and (\ref{u}) imply that the
Christoffel tensors $\lambda$ and $\tilde\lambda$ associated to
these adapted frames are equal at corresponding points.
We then conclude that $A$ and $\tilde A$ both satisfy the equation
\begin{equation}
A^{-1}\dd A = \omega-\lambda
\end{equation}
Now, left translation by $f(x_0)\tilde f(x_0)^{-1}$ followed by a
suitable rotation in $T_{f(x_0)}N$, if necessary, assure that we may
suppose that $A(x_0)=\tilde A(x_0)$. Hence, the uniqueness of
Darboux primitives in a simply connected domain implies that
$A=\tilde A$.

Thus, we have
\begin{equation}
\omega_{\mathfrak{m}}|_{f(x)}(f_* e_a) = \sum_k \bar e_k A^k_a
\end{equation}
and
\begin{equation}
\omega_{\mathfrak{m}}|_{\tilde f(x)}(\tilde f_* e_a) =\sum_k \bar
e_k A^k_a.
\end{equation}
Therefore, $f$ and $\tilde f$ describe integral leaves of the
distribution $\mathcal{P}$ we defined above passing through the
point $f(x_0)\in N$. The uniqueness part of Frobenius's theorem
implies that $f=\tilde f$.

This finishes the proof of Theorem \ref{main}.


\vspace{2cm}
\noindent Jorge H. S. de Lira \\
(corresponding author)\\
Departamento de Matem\'atica - UFC\\
Campus do Pici, Bloco 914\\
Fortaleza, Cear\'a, Brazil\\
60455-760\\
jorge.lira@pq.cnpq.br

\vspace{0.5cm}

\noindent Marcos F. de Melo\\
UFC - Campus do Cariri\\
Av. Ten. Raimundo Rocha\\
Juazeiro de Norte, Cear\'a, Brazil\\
60030-200\\
mcosmelo79@hotmail.com\\

\end{document}